\documentclass[12pt, reqno]{amsart}
\usepackage{amsmath, amstext, amsbsy, amssymb, amscd, graphics}

\newtheorem{theorem}{Theorem}[section]
\newtheorem{lemma}[theorem]{Lemma}
\newtheorem{proposition}[theorem]{Proposition}

\theoremstyle{definition}
\newtheorem{definition}[theorem]{Definition}

\theoremstyle{remark}
\newtheorem{remark}[theorem]{Remark}

\numberwithin{equation}{section}

\newcommand{\C}{ \mathbb C }

\newcommand{\End}{{\rm End}}
\newcommand{\fock}{{\mathbb H}_X}
\newcommand{\g}{{\gamma}}
\newcommand{\Hn}{H^*(\Xn)}
\newcommand{\Hx}{{\mathbb H}_X}

\newcommand{\la}{\lambda}

\newcommand{\Q}{ \mathbb Q }
\newcommand{\Supp}{{\rm Supp}}

\newcommand{\vac}{|0\rangle}
\newcommand{\w}{\tilde}
\newcommand{\W}{\widetilde}

\newcommand{\Xn}{ X^{[n]}}
\newcommand{\Z}{ \mathbb Z }

{\vskip-\lastskip\medskip
  \noindent
  {\em #1.}\enspace
  }%
{\qed\par\medskip
  }

\begin{document}
\title[Integral operators and integral classes]
{Integral operators and integral cohomology classes of
Hilbert schemes}

\author{Zhenbo Qin${}^\dagger$}
\address{Department of Mathematics, University of Missouri,
Columbia, MO 65211, USA} \email{zq@math.missouri.edu}
\thanks{${}^\dagger$Partially supported by an NSF grant}

\author{Weiqiang Wang}
\address{Department of Mathematics, University of Virginia,
Charlottesville, VA 22904} \email{ww9c@virginia.edu}


\begin{abstract}
The methods of integral operators on the cohomology of Hilbert
schemes of points on surfaces are developed. They are used to
establish integral bases for the cohomology groups of Hilbert
schemes of points on a class of surfaces (and conjecturally, for
all simply connected surfaces).
\end{abstract}

\maketitle
\date{}

\section{Introduction}

The Hilbert scheme $\Xn$ parameterizes all length-$n$
$0$-dimensional closed subscheme of a complex smooth projective
surface $X$. A classical result of G\"ottsche \cite{Got}
calculated the Betti numbers of the Hilbert scheme $\Xn$ for
an arbitrary surface $X$.
Nakajima \cite{Na1} (see also \cite{Gro}) constructed a
Heisenberg algebra which acts irreducibly on the direct sum $\Hx$
of the rational cohomology of the Hilbert schemes $\Xn$ for all
$n$. As a corollary, a linear basis for the {\it rational}
cohomology of $\Xn$ in terms of the Heisenberg algebra operators
can be constructed. With the help of this fundamental
construction, there has been intensive activities and significant
progress in the study of Hilbert schemes such as their rational
cohomology rings made by Lehn, Sorger,
W.-P. Li and the authors of this
paper, and others (cf. \cite{QW} for extensive references).
However, there has been virtually no work toward the basic
problem of studying the {\em integral} cohomology of $\Xn$ for
a general surface $X$ (on the other hand, see \cite{ES2, LS} when
$X$ is the projective or affine plane).

The purpose of this paper is to develop some general techniques of
integral operators and then to use them to find integral bases for
the {\it integral} cohomology of $\Xn$ for a certain class of
surfaces $X$. By an integral operator we mean a linear operator on
$\Hx$ which sends every integral cohomology class to an integral
one. One of the difficulties of studying the integral cohomology
of $\Xn$ is that not many rational cohomology classes are known to
be integral. The starting point of this work is the following key
observation which we derive from the Stability Theorem established
in \cite{LQW}: if $A$ is an integral cohomology class of $\Xn$ and
if $A$ is written as $\mathfrak a_A \vac$ where $\mathfrak a_A$ is
a polynomial of creation Heisenberg operators, then $\mathfrak
a_A$ is an integral operator (Proposition~\ref{prop_int}). This
provides an effective method of constructing new integral classes
from known ones.

Our study of integral operators and integral classes are roughly
divided into two distinct parts. The first part involves integral
operators and integral classes that are `created' from the
creation Heisenberg operators $\mathfrak a_{-r}(1)$ and $\mathfrak
a_{-r}(x)$ associated to the identity cohomology class $1$
and the point cohomology class $x$.
For a class $\alpha$ in the rational cohomology
$H^*(X)$ and for a partition $\la =(1^{m_1}2^{m_2} \ldots)$ where
$m_r$ stands for the number of parts equal to $r$, define
$|\la| = \sum_{r \ge 1} r m_r$, $\ell(\la) = \sum_{r \ge 1} m_r$,
\begin{eqnarray}
{\mathfrak z}_{\la}
  = \prod_{r \ge 1}  r^{m_r} m_r!, \quad
\mathfrak a_{-\la}(\alpha) = \prod_{r \ge 1}
  {\mathfrak a_{-r}(\alpha)^{m_r}}.     \label{intro_for}
\end{eqnarray}
We prove that the operator $1/{\mathfrak z}_{\la} \cdot \mathfrak
a_{-\la}(1)$ is integral. By construction, the operator $\mathfrak
a_{-\la}(x)$ associated to the point class $x$ is always integral.

The second part of our study is to develop integral operators and
integral classes created from the creation Heisenberg operators
$\mathfrak a_{-i}(\alpha)$ associated with $\alpha \in H^2(X)$.
Recall from \cite{Na2} (see also \cite{Gro}) the subvarieties
$L^\la C \subset \Xn$ associated to partitions $\la$ of $n$ and an
embedded curve $C$ (see Subsection~\ref{subsect_c} for definitions).
Let $\mathfrak m_{\la, C}$ be the integral operator $\mathfrak
a_{[L^\la C]}$. We first extend the definitions of $[L^\la C]$ and
$\mathfrak m_{\la, C}$ to $[L^\la \alpha]$ and $\mathfrak m_{\la,
\alpha}$ for an arbitrary class $\alpha \in H^2(X)$. We prove that
$\mathfrak m_{\la, \alpha}$ is an integral operator for $\alpha
=\alpha_1 \pm \alpha_2$ whenever $\mathfrak m_{\la, \alpha_1}$ and
$\mathfrak m_{\la, \alpha_2}$ are integral. In particular, it
follows that $\mathfrak m_{\la, \alpha}$ is integral for every
divisor $\alpha$ (Theorem~\ref{div}). Moreover, we show
(cf.~Theorem~\ref{unimod}) that if the intersection matrix of
$\alpha_1, \cdots, \alpha_k$ in $H^2(X)$ has determinant $\pm 1$,
then so is the determinant of the intersection matrix of the
classes $ \mathfrak m_{\nu^1, \alpha_1} \cdots \mathfrak m_{\nu^k,
\alpha_k} \vac$ in $H^*(\Xn)$, where the partitions $\nu^i$
satisfy $|\nu^1| + \cdots + |\nu^k| = n$. Our proofs of these
results use in an effective way the interrelations among the
operators $\mathfrak m_{\nu, \alpha}$,
the creation Heisenberg operators,
and the monomial symmetric functions $m_\nu$ in
the ring of symmetric functions (cf. Chapter~9, \cite{Na2}).

As an application of the above development, we obtain the
following.

\begin{theorem} \label{intro_thm}
Let $X$ be a projective surface with $H^1(X; \mathcal O_X) =
H^2(X; \mathcal O_X) = 0$. Let $\alpha_1, \cdots, \alpha_k$ be an
integral basis for $H^2(X; \mathbb Z)/\text{\rm Tor}$.
Then the following classes
\begin{eqnarray*}
\frac{1}{{\mathfrak z}_\la} \cdot \mathfrak a_{-\la}(1) \mathfrak
a_{-\mu}(x) \mathfrak m_{\nu^1, \alpha_1} \cdots \mathfrak
m_{\nu^k, \alpha_k} \vac, \quad |\la|+|\mu|+ \sum_{i=1}^k |\nu^i|
= n
\end{eqnarray*}
are integral, and furthermore, they form an integral basis for
$H^*(\Xn; \mathbb Z)/\text{\rm Tor}$.
\end{theorem}

We conjecture that the cohomology class $[L^\la \alpha]$ is
integral whenever $\alpha \in H^2(X)$ is an integral class.
If this conjecture is true, then the statement in
Theorem~\ref{intro_thm} will be valid for
every simply connected projective surface $X$.
We refer to Remark~\ref{rmk_monodromy} for some discussions
in this direction.

In view of the dictionary between Hilbert schemes and symmetric
products developed in \cite{QW}, the counterparts of the results
in this paper are expected to be valid for the Chen-Ruan orbifold
cohomology rings of the symmetric products.

Our paper is organized as follows. In Sect.~\ref{sect_basics}, we
review some basic results about the Hilbert schemes. In
Sect.~\ref{sect_int} and Sect.~\ref{sect_deg_2}, integral
operators related to the creation
operators $\mathfrak a_{-i}(1)$ and $\mathfrak a_{-i}(\alpha)$
with $\alpha \in H^2(X)$ are investigated respectively.
In Sect.~\ref{sect_app}, we prove Theorem~\ref{intro_thm}
and some other general structure results.

\medskip\noindent
{\bf Conventions.} Unless otherwise indicated, all the cohomology
groups in this paper are in $\mathbb Q$-coefficients. For a
continuous map $p: Y_1 \to Y_2$ between two smooth compact
manifolds and for $\alpha_1 \in H^*(Y_1)$, we define
$p_*(\alpha_1)$ to be ${\rm PD}^{-1}p_{*}({\rm PD}(\alpha_1))$
where ${\rm PD}$ stands for the Poincar\'e duality. For a smooth
projective surface $X$, by abusing notations, we let $1 \in H^0(X;
\Z)$ be the fundamental cohomology class of $X$; in addition, we
let $x$ denote either a point in $X$, or the class in $H^4(X; \Z)$
which is the Poincar\'e dual of the homology class represented by
a point in $X$.

\bigskip\noindent
{\bf Acknowledgments.} The authors thank Wei-Ping Li for valuable
discussions, and are very grateful to the referee for insightful
comments and suggestions which lead to
Remark~\ref{rmk_monodromy} and subsection~\ref{sec:model}.
After communicating this paper to E. Markman, we learned about his
preprint \cite{Mar} in which among other interesting results the
{\it ring generators} for the integral cohomology rings of Hilbert
schemes of points on $K3$ surfaces are obtained.

\section{\bf Basics on Hilbert schemes of points on surfaces}
\label{sect_basics}

Let $X$ be a complex smooth projective surface,
and $\Xn$ be the Hilbert scheme of points in $X$.
An element in $\Xn$ is represented by a
length-$n$ $0$-dimensional closed subscheme $\xi$ of $X$. For $\xi
\in \Xn$, let $I_{\xi}$ be the corresponding sheaf of ideals. It
is well known that $\Xn$ is smooth. Sending an element in $\Xn$ to
its support in the symmetric product ${\rm Sym}^n(X)$, we obtain
the Hilbert-Chow morphism $\pi_n: \Xn \rightarrow {\rm Sym}^n(X)$,
which is a resolution of singularities. Let ${\mathcal Z}_n$
be the universal codimension-$2$ subscheme of $\Xn\times X$,
which can be described set-theoretically by
\begin{eqnarray}  \label{cod2}
{\mathcal Z}_n=\{(\xi, x) \subset \Xn\times X \, | \, x\in
{\rm Supp}{(\xi)}\}\subset \Xn\times X.
\end{eqnarray}
Let $\Hn$ be the total cohomology of $\Xn$ with $\Q$-coefficients,
and put
\begin{eqnarray}  \label{fock}
\fock = \bigoplus_{n=0}^\infty \Hn.
\end{eqnarray}

For $m \ge 0$ and $n > 0$, let $Q^{[m,m]} = \emptyset$,
and let $Q^{[m+n,m]}$ be the closed subscheme of
$X^{[m+n]} \times X \times X^{[m]}$ defined in \cite{Lehn}
whose set-theoretical description is:
\begin{eqnarray*}
\{ (\xi, x, \eta) \in X^{[m+n]} \times X \times X^{[m]} \, | \,
\xi \supset \eta \text{ and } \mbox{Supp}(I_\eta/I_\xi)
= \{ x \} \}.
\end{eqnarray*}

We recall Nakajima's definition of the Heisenberg operators
\cite{Na1}. Let $n > 0$. The linear operator $\mathfrak
a_{-n}(\alpha) \in \End(\fock)$ with $\alpha \in H^*(X)$ is
defined by
\begin{eqnarray}  \label{def_a}
\mathfrak a_{-n}(\alpha)(A) = \tilde{p}_{1*}([Q^{[m+n,m]}] \cdot
\tilde{\rho}^*\alpha \cdot \tilde{p}_2^*A)
\end{eqnarray}
for $A \in H^*(X^{[m]})$, where $\tilde{p}_1, \tilde{\rho},
\tilde{p}_2$ are the projections of $X^{[m+n]} \times X \times
X^{[m]}$ to $X^{[m+n]}, X, X^{[m]}$ respectively. Define
$\mathfrak a_{n}(\alpha) \in \End(\fock)$ to be $(-1)^n$ times the
operator obtained from the definition of $\mathfrak
a_{-n}(\alpha)$ by switching the roles of $\tilde{p}_1$ and $
\tilde{p}_2$. We often refer to $\mathfrak a_{-n}(\alpha)$ (resp.
$\mathfrak a_n(\alpha)$) as the {\em creation} (resp. {\em
annihilation})~operator. We also set $\mathfrak a_0(\alpha) =0$. A
non-degenerate super-symmetric bilinear form $(-, -)$ on $\mathbb
H_X$ is induced from the standard one on $H^*(\Xn)$ defined by
\begin{eqnarray*}
(\alpha,\beta) =\int_{\Xn} \alpha\beta, \qquad
\alpha, \beta\in H^*(\Xn).
\end{eqnarray*}
This allows us to define the {\it adjoint} $\mathfrak f^\dagger
\in \End(\mathbb H_X)$ for $\mathfrak f \in \End(\mathbb H_X)$.
Then,
\begin{eqnarray} \label{a_n_pm}
\mathfrak a_n(\alpha) = (-1)^n \cdot \mathfrak
a_{-n}(\alpha)^\dagger.
\end{eqnarray}

The operators $\mathfrak a_{n}(\alpha) \in \End(\fock)$ with $\alpha
\in H^*(X)$ and $n \in \mathbb Z$ satisfy the following
Heisenberg algebra commutation relation (cf. \cite{Na2}):
\begin{eqnarray}  \label{eq:heis}
[\mathfrak a_m(\alpha), \mathfrak a_n(\beta)] = -m \cdot
\delta_{m,-n} \cdot (\alpha, \beta) \cdot {\rm Id}_{\fock}.
\end{eqnarray}
The space $\fock$ is an irreducible
module over the Heisenberg algebra generated by the operators
$\mathfrak a_n(\alpha)$ with a highest~weight~vector
\begin{eqnarray*}
\vac=1 \in H^0(X^{[0]}) \cong \mathbb Q.
\end{eqnarray*}
It follows that $\fock$ is linearly
spanned by all the {\it Heisenberg monomial classes}:
\begin{eqnarray}  \label{Heis_mon}
\mathfrak a_{-n_1}(\alpha_1) \cdots
\mathfrak a_{-n_k}(\alpha_k)\vac
\end{eqnarray}
where $k \ge 0, n_1, \ldots, n_k > 0$, and $\alpha_1,
\ldots, \alpha_k$ run over a linear basis of $H^*(X)$.

For a nonnegative integer $n$, we define the operator:
\begin{eqnarray}
   {\bf 1}_{-n}
=\frac{1}{n!} \cdot \mathfrak a_{-(1^n)}(1)
   = \frac{1}{n!} \cdot \mathfrak a_{-1}(1)^n. \label{a_la_3}
\end{eqnarray}
The geometric meaning of ${\bf 1}_{-n}$ is that
${\bf 1}_{-n} \vac$ is equal to the fundamental class of
the Hilbert scheme $X^{[n]}$.
For simplicity, we set ${\bf 1}_{-n} = 0$ when $n < 0$.

The following is the Stability Theorem 5.1 proved in \cite{LQW}.

\begin{theorem} \label{cup_product}
Let $s \ge 1$ and $k_i \ge 1$ for $1 \le i \le s$. Fix $n_{i, j}
\ge 1$ and $\alpha_{i, j} \in H^*(X)$ for $1 \le j \le k_i$, and
fix $n$ with $n \ge \sum\limits_{j=1}^{k_i} n_{i, j}$ for all $1 \le i
\le s$. Then the cup product
\begin{eqnarray} \label{cup_product1}
\prod_{i=1}^s \left ( {\bf 1}_{-(n - \sum_{j=1}^{k_i} n_{i, j})}
\left (\prod_{j=1}^{k_i} \mathfrak a_{-n_{i, j}}(\alpha_{i, j})
\right ) |0\rangle \right )
\end{eqnarray}
in $H^*(\Xn)$ is equal to a finite linear combination of
monomials of the form
\begin{eqnarray} \label{cup_product2}
 {\bf 1}_{-(n - \sum_{p=1}^N m_{p})}
\left ( \prod_{p=1}^N \mathfrak a_{-m_{p}}(\g_{p}) \right )
|0\rangle
\end{eqnarray}
whose coefficients are independent of $X,
\alpha_{i,j}$ and $n$. Here $\sum\limits_{p=1}^N m_{p} \le
\sum\limits_{i=1}^s \sum\limits_{j=1}^{k_i}
n_{i,j}$, and $\g_1, \ldots, \g_N$ depend only on $\alpha_{i,j}$,
$1\le i \le s, 1\le j\le k_i$, and the canonical class and
the Euler class of $X$.  In addition,
the expression (\ref{cup_product2}) satisfies the upper bound
$$\sum\limits_{p=1}^N m_{p} = \sum\limits_{i=1}^s
\sum\limits_{j=1}^{k_i} n_{i,j}$$ if and only if it is
$\displaystyle{{\bf 1}_{-(n -\sum_{i=1}^s \sum_{j=1}^{k_i}
n_{i,j})} \left ( \prod_{i=1}^s\prod_{j=1}^{k_i} \mathfrak
a_{-n_{i, j}}(\alpha_{i, j}) \right ) |0\rangle}$ with coefficient
$1$.
\end{theorem}

\section{\bf Integral operators involving only $1 \in H^*(X)$}
\label{sect_int}
\subsection{\bf Integral operators}
\label{subsect_int}

\begin{definition} \label{def_int}
(i) A class $A \in H^*(\Xn)$ is {\it integral} if it is
contained in
\begin{eqnarray*}
H^*(\Xn; \Z)/\text{Tor} \subset H^*(\Xn);
\end{eqnarray*}

(ii) A linear basis of $H^*(\Xn)$ is {\it integral} if
its members are integral classes and form a $\Z$-basis of
the lattice $H^*(\Xn; \Z)/\text{Tor}$;

(iii) A linear operator $\mathfrak f \in \End(\mathbb H_X)$ is
{\it integral} if $\mathfrak f(A) \in \mathbb H_X$ is an integral
class whenever $A \in \mathbb H_X$ is an integral cohomology
class.
\end{definition}

A linear basis of $H^*(\Xn)$ is integral if
and only if its members are integral classes and the matrix
formed by the pairings of its members is unimodular.

\begin{lemma} \label{int_Heis}
{\rm (i)} If $\mathfrak f \in \End(\mathbb H_X)$ is integral, then
so is its adjoint $\mathfrak f^\dagger$;

{\rm (ii)} The Heisenberg operators $\mathfrak a_{n}(\alpha), n
\in \mathbb Z$ are integral if $\alpha \in H^*(X)$ is integral.
\end{lemma}
\begin{proof}
(i) Note that a class $A \in \mathbb H_X$ is integral if and only
if $(A, B)$ is an integer whenever $B \in \mathbb H_X$ is an
integral class. It follows that $\mathfrak f \in \End(\mathbb
H_X)$ is integral if and only if its adjoint operator $\mathfrak
f^\dagger \in \End(\mathbb H_X)$ is integral.

(ii) Recall that $\mathfrak a_0(\alpha) = 0$. Next, fix $n > 0$.
By (\ref{def_a}), the Heisenberg operator $\mathfrak
a_{-n}(\alpha)$ is integral. By (\ref{a_n_pm}) and (i), the
operator $\mathfrak a_{n}(\alpha)$ is integral as well.
\end{proof}

\begin{lemma} \label{int_1_n}
For $n \ge 0$, the operator ${\bf 1}_{-n}$ is integral.
\end{lemma}
\begin{proof}
Recall from (\ref{a_la_3}) that ${\bf 1}_{-n}
= 1/n! \cdot \mathfrak a_{-1}(1)^n$. Fix any integer $m \ge 0$
and an integral class $A \in H^j(X^{[m]})$.
For $i = m, \ldots, m+n -1$, let $Q_i$ be the image of
the subscheme $Q^{[i+1, i]}$ under the natural projection
$X^{[i+1]} \times X \times X^{[i]} \to X^{[i+1]} \times X^{[i]}$.
Set-theoretically, $Q_i = \{ (\xi_{m+1}, \xi_m) \in X^{[i+1]}
\times X^{[i]} \, | \, \xi_{m+1} \supset \xi_m \}$.
Let $\phi_{i, 1}$ and $\phi_{i, 2}$ be the two projections
of $X^{[i+1]} \times X^{[i]}$. It follows from (\ref{def_a}) that
\begin{eqnarray*}
\mathfrak a_{-1}(1)(A) = (\phi_{m, 1})_*([Q_m] \cdot \phi_{m, 2}^*A).
\end{eqnarray*}
Repeating this process and using the projection formula,
we conclude that
\begin{eqnarray*}
\mathfrak a_{-1}(1)^n(A) = (\w \phi_{m+n})_*([Q] \cdot \w \phi_{m}^*A)
\end{eqnarray*}
where $\w \phi_{i}$ denotes the projection of
$Y := X^{[m+n]} \times \cdots \times X^{[m+1]} \times X^{[m]}$
to $X^{[i]}$ and
\begin{eqnarray*}
Q = \{(\xi_{m+n}, \ldots, \xi_{m+1}, \xi_m) \in Y|\,\,
\xi_{m+n} \supset \cdots \supset \xi_{m+1} \supset \xi_m\}
\end{eqnarray*}
(the scheme structure on $Q$ can be described similarly as
that on $Q_i$).
To show that ${\bf 1}_{-n}(A) \in H^j(X^{[m+n]})$ is integral,
it suffices to prove that the intersection number
$(\w \phi_{m+n})_*([Q] \cdot \w \phi_{m}^*A) \cdot B$
is divisible by $n!$ for any integral class
$B \in H^{4(m+n)-j}(X^{[m+n]})$.

Represent the integral class $A$ by a piecewise smooth cycle
$W_A \subset X^{[m]}$. Then, the integral class
$[Q] \cdot \w \phi_{m}^*A$ is represented by
\begin{eqnarray*}
Q_A := \{(\xi_{m+n}, \ldots, \xi_{m+1}, \xi_m) \in Y|\,\,
\xi_{m+n} \supset \cdots \supset \xi_{m+1} \supset \xi_m
\text{ and } \xi_m \in W_A\}.
\end{eqnarray*}
Note that $\w \phi_{m+n}|_{Q_A}: Q_A \to \w \phi_{m+n}(Q_A)$
is generically finite, and a generic element in
$\w \phi_{m+n}(Q_A)$ is of the form $\xi_m + x_1+\ldots + x_n$
where $\xi_m \in W_A$ is generic, the points $x_1, \ldots, x_n$
are distinct, and $\text{Supp}(\xi_m) \cap \{x_1, \ldots, x_n\}
= \emptyset$.

Represent $B$ by a piecewise smooth cycle
$W_B$ such that $W_B$ intersects $\w \phi_{m+n}(Q_A)$
transversely at generic points $P_1, \ldots, P_s$.
Write $P_i$ as $\xi_m + x_1+\ldots + x_n$. Let
\begin{eqnarray*}
\xi_m = \eta_1 + \ldots + \eta_t + x_1' + \ldots + x_u'
+x_1''+\ldots +x_v''
\end{eqnarray*}
where each $\eta_i$ is supported at one point with
$\ell(\eta_i) \ge 2$, the following subsets of $X$:
\begin{eqnarray*}
\Supp(\eta_1), \ldots, \Supp(\eta_t), \{x_1'\}, \ldots,
\{x_u'\},  \{x_1''\}, \ldots, \{x_v''\}
\end{eqnarray*}
are mutually disjoint, and only the points $x_1', \ldots, x_u'$
can move in Zariski open subsets of the surface $X$.
(In other words, if we put
\begin{eqnarray*}
U = X - \Supp(\eta_1) \cup \ldots \cup \Supp(\eta_t)
\cup \{x_1', \ldots, x_u', x_1'', \ldots, x_v''\},
\end{eqnarray*}
then for all distinct points $\w x_1', \ldots, \w x_u'$ in $U$,
we have
\begin{eqnarray*}
\eta_1 + \ldots + \eta_t + \w x_1' + \ldots + \w x_u'
+x_1''+\ldots +x_v'' \in W_A. )
\end{eqnarray*}
So $\big ( \w \phi_{m+n}|_{Q_A} \big )^{-1}(P_i)$ consists
of $(u+n)(u+n-1)\cdots (u+1)$ distinct points.
It follows that the intersection number
$(\w \phi_{m+n})_*([Q] \cdot \w \phi_{m}^*A) \cdot B$
is divisible by $n!$.
\end{proof}

Let $A \in H^*(\Xn)$. By (\ref{Heis_mon}), $A$ is
a linear combination of classes of the form
$\mathfrak a_{-1}(1)^m \mathfrak a_{-m_1}(\alpha_1)\cdots
\mathfrak a_{-m_\ell}(\alpha_\ell)\vac$
where ${\alpha_1, \ldots, \alpha_\ell \in
\bigoplus_{i \ge 1} H^i(X)}$, $m \ge 0$, $\ell \ge 0$,
and $m_1, \ldots, m_\ell > 0$.
By (\ref{eq:heis}),
$\mathfrak a_1(x) \big (\mathfrak a_{-m_1}(\alpha_1)\cdots
\mathfrak a_{-m_\ell}(\alpha_\ell)\vac \big ) = 0$.
It follows that the class $A \in H^*(\Xn)$ can be written as
\begin{eqnarray} \label{Ai}
A = {\bf 1}_{-n_1}(A_1) + \ldots + {\bf 1}_{-n_k}(A_k)
\end{eqnarray}
where $0 \le n_1 < \ldots < n_k$, and $A_i \in H^*(X^{[n-n_i]})$
with $\mathfrak a_1(x)(A_i) = 0$ for every $i$.

\begin{lemma} \label{int_Ai}
Let $A \in H^*(\Xn)$ be expressed as in (\ref{Ai}). Then,
\begin{enumerate}
\item[{\rm (i)}] the classes $A_i$ are uniquely
determined by $A$;

\item[{\rm (ii)}] $A$ is integral if and only if all
the classes $A_i$ are integral.
\end{enumerate}
\end{lemma}
\begin{proof}
First of all, by Lemma~\ref{int_1_n}, $A$ is integral if
all the classes $A_i$ are integral. So it remains to
prove (i) and the ``only if" part of (ii).

Next, applying the operator $\mathfrak a_1(x)^{n_k}$ to both sides
of (\ref{Ai}), we obtain
\begin{eqnarray*}
\mathfrak a_1(x)^{n_k}(A) = \mathfrak a_1(x)^{n_k} \big ({\bf
1}_{-n_1}(A_1) + \ldots + {\bf 1}_{-n_k}(A_k) \big )
= (-1)^{n_k} \cdot A_k.
\end{eqnarray*}
Thus the class $A_k$ is uniquely determined by $A$.
Moreover, by Lemma~\ref{int_Heis}~(ii), if $A$ is
an integral class, then so is the class $A_k$.

Finally, repeating the above process to the class
$A' := (A - {\bf 1}_{-n_k}(A_k))$, we conclude that
(i) and the ``only if" part of (ii) hold for
all the classes $A_i$.
\end{proof}

\begin{proposition} \label{prop_int}
Let $A \in \fock$ be an integral class. Write $A$ as $\mathfrak
a_A \vac$ where $\mathfrak a_A$ is a unique polynomial of creation
operators. Then, $\mathfrak a_A$ is an integral operator.
\end{proposition}
\begin{proof}
Let $A \in H^*(\Xn)$. Fix an integral class $B \in H^*(X^{[m]})$.
We want to show that the cohomology class $\mathfrak a_A(B)
\in H^*(X^{[m+n]})$ is still integral.

First of all, decomposing $A$ as in (\ref{Ai}), we see from
Lemma~\ref{int_Ai}~(ii) and Lemma~\ref{int_1_n} that we may assume
$\mathfrak a_1(x)(A) = 0$. Similarly, write $B
= {\bf 1}_{-m_1}(B_1) + \ldots + {\bf 1}_{-m_k}(B_k)$ as in (\ref{Ai}),
where $0 \le m_1 < \ldots < m_k \le m$, and $B_i \in H^*(X^{[m-m_i]})$
with $\mathfrak a_1(x)(B_i) = 0$ for every $i$.
By Lemma~\ref{int_Ai}~(ii), each class $B_i$ is integral. Now
\begin{eqnarray*}
  \mathfrak a_A(B)
= \sum_{i=1}^k \mathfrak a_A \big ({\bf 1}_{-m_i}(B_i) \big )
= \sum_{i=1}^k {\bf 1}_{-m_i} \big (\mathfrak a_A (B_i) \big ).
\end{eqnarray*}
It follows from Lemma~\ref{int_1_n} that we may assume
$\mathfrak a_1(x)(B) = 0$ as well.

By Lemma~\ref{int_1_n} again, we have two
integral classes ${\bf 1}_{-m}A, {\bf 1}_{-n}B \in H^*(X^{[m+n]})$.
By Theorem~\ref{cup_product}, the cup product
$({\bf 1}_{-m}A) \cdot ({\bf 1}_{-n}B)$ equals
\begin{eqnarray*}
\mathfrak a_A(B) + \sum_{i
\ge 1} {\bf 1}_{-i}(A_i)
\end{eqnarray*}
where $\mathfrak a_1(x)(A_i) = 0$ for every $i \ge 1$. By our
assumption on $A$ and $B$, we have $\mathfrak a_1(x) \big
(\mathfrak a_A(B) \big ) = 0$. Hence $\mathfrak a_A(B) \in
H^*(X^{[m+n]})$ is integral by Lemma~\ref{int_Ai}.
\end{proof}
\subsection{\bf Integral classes and operators involving
only $1 \in H^*(X)$}
\label{subsect_1}
\par
$\,$

In view of Proposition~\ref{prop_int}, to obtain integral
operators, we need to construct integral classes.
We begin with the Chern classes of some tautological
vector bundles over $X^{[n]}$. For a line bundle $L$ on
the surface $X$, let $L^{[n]} = \big (p_1|_{{\mathcal Z}_n}\big )_*
\big (p_2|_{{\mathcal Z}_n}\big )^*L$
where ${\mathcal Z}_n$ is from (\ref{cod2}),
and $p_1$ and $p_2$ are the projections of $\Xn \times X$ to
its two factors. By the Theorem~4.6 in \cite{Lehn}
which was conjectured earlier by G\" ottsche,
\begin{eqnarray} \label{c_i}
   c_i(\mathcal O_X^{[n]})
= (-1)^i \cdot \quad \sum_{|\nu|=n, \,\, \ell(\nu) = n-i}
  \frac{\mathfrak a_{-\nu}(1)\vac}{{\mathfrak z}_\nu}.
\end{eqnarray}

\begin{lemma} \label{int_a_la}
For every partition $\la$, we have

{\rm (i)} $1/{\mathfrak z}_\la \cdot \mathfrak a_{-\la}(1)
\vac$ is an integral class;

{\rm (ii)} $1/{\mathfrak z}_\la \cdot \mathfrak a_{-\la}(1)$
is an integral operator.
\end{lemma}
\begin{proof}
(ii) follows from (i) and Proposition~\ref{prop_int}.
To prove (i), we let $n = |\la|$ and use induction on $n$.
Our result is trivially true when $n = 0, 1$.
In the following, assume that
$1/{\mathfrak z}_\mu \cdot \mathfrak a_{-\mu}(1) \vac$ is
integral whenever $|\mu| < n$.

Let $\la = (1^{m_1}2^{m_2}3^{m_3} \cdots)$ so that $n
= \sum_r rm_r$. First of all, assume $m_r > 0$ for at least
two different $r$'s. Then, $rm_r < n$ for every $r$.
Putting $\mu^r = (r^{m_r})$ for $r \ge 1$ and
applying induction to the partitions $\mu^r$, we obtain
integral classes $A_r := 1/{\mathfrak z}_{\mu^r} \cdot
\mathfrak a_{-\mu^r}(1) \vac$. By Proposition~\ref{prop_int},
the operators $\mathfrak a_{A_r} = 1/{\mathfrak z}_{\mu^r}
\cdot \mathfrak a_{-\mu^r}(1)$ are integral. Thus,
$1/{\mathfrak z}_\la \cdot \mathfrak a_{-\la}(1) \vac =
\prod_{r \ge 1} \big ( 1/{\mathfrak z}_{\mu^r}
\cdot \mathfrak a_{-\mu^r}(1) \big ) \cdot \vac$ is integral.

We are left with the case when $m_r > 0$ for a unique $r$.
In this case, $|\la| = rm_r = n$ and $\ell(\la) = m_r$.
Applying (\ref{c_i}) to $i = n - m_r$, we have
\begin{eqnarray} \label{int_a_la1}
   c_i(\mathcal O_X^{[n]})
&=&(-1)^i \cdot \quad \sum_{|\nu|=rm_r, \,\, \ell(\nu) = m_r}
   \frac{\mathfrak a_{-\nu}(1)\vac}{{\mathfrak z}_\nu}
   \nonumber  \\
&=&(-1)^i \cdot \frac{\mathfrak a_{-\la}(1) \vac}
   {{\mathfrak z}_\la} + (-1)^i \cdot
   \sum_{\stackrel{\nu \ne \la, |\nu|=rm_r}{\ell(\nu) = m_r}}
   \frac{\mathfrak a_{-\nu}(1)\vac}{{\mathfrak z}_\nu}.
\end{eqnarray}
Note that for any partition $\nu = (1^{t_1}2^{t_2} \cdots)$
satisfying $\nu \ne \la, |\nu|= rm_r$ and $\ell(\nu) = m_r$,
there are at least two $i$'s with $t_i > 0$.
By the previous paragraph,
$1/{\mathfrak z}_\nu \cdot \mathfrak a_{-\nu}(1) \vac$ is
integral. Hence we see from (\ref{int_a_la1}) that
$1/{\mathfrak z}_\la \cdot \mathfrak a_{-\la}(1) \vac$
is integral.
\end{proof}

\section{\bf Integral operators involving only classes in
$H^2(X)$}
\label{sect_deg_2}

\subsection{\bf Integral classes involving only an embedded curve
in $X$}
\label{subsect_c}
\par
$\,$

Let $C$ be a smooth irreducible curve in the surface $X$.
By abusing notations,
we also use $C$ to denote the corresponding
divisor and the corresponding cohomology class.
Following Subsection 9.3 of \cite{Na2}
(see also \cite{Gro}), for every partition
$\la = (\la_1, \la_2, \ldots )$ of $n$, we define the subvariety
$L^\la C := \overline{(\pi_n)^{-1}(S^n_\la C)}$ of $\Xn$,
where $S^n_\la C = \{\sum_i \la_i x_i| \,\, x_i \in C,
x_i \ne x_j \text{ for } i \ne j\}$,
and $\pi_n$ is the Hilbert-Chow morphism.
For $n \ge 0$, let $\mathbb H_{n, C}$ be the $\mathbb Q$-linear
span of all the classes $\mathfrak a_{-\la}(C)\vac, \la \vdash n$
where $\la \vdash n$ denotes that $\la$ is a partition of $n$.
By the Theorem 9.14 in \cite{Na2}, the integral class
$[L^\la C] \in H^{2n}(X^{[n]}; \Z)$ is contained in
$\mathbb H_{n, C} \subset H^*(\Xn)$. Define
\begin{eqnarray} \label{H_C}
\mathbb H_{C} = \bigoplus_{n=0}^\infty \mathbb H_{n, C}.
\end{eqnarray}

Let $\Lambda$ be the ring of symmetric functions in infinitely
many variables (see p.19 of \cite{Mac}),
and $\Lambda_\Q = \Lambda \otimes_\Z \Q$. Let $\Lambda^n$ and
$\Lambda_\Q^n$ be the degree-$n$ parts in $\Lambda$ and
$\Lambda_\Q$ respectively. For a partition $\la$,
let $p_\la, m_\la$ and $s_\la$ be the power-sum symmetric
function, the monomial symmetric function and
the Schur function respectively.

In \cite{Na2}, Nakajima defined a linear isomorphism
\begin{eqnarray} \label{Phi_C1}
\Phi_C: \quad \Lambda_\Q = \bigoplus_{n=0}^\infty \Lambda_\Q^n
\quad \to \quad \mathbb H_{C}
\end{eqnarray}
which satisfies the following two properties:
\begin{eqnarray} \label{Phi_C2}
\Phi_C(p_\la) = \mathfrak a_{-\la}(C)\vac, \qquad
\Phi_C(m_\la) = [L^\la C].
\end{eqnarray}

For a partition $\la$, define $\mathfrak m_{\la, C} =
\mathfrak a_{[L^\la C]} \in \End(\mathbb H_X)$.
By Proposition~\ref{prop_int}, the operator $\mathfrak m_{\la, C}$
is integral. Moreover, we see from (\ref{Phi_C2}) that
$\mathfrak m_{\la, C}$ is a polynomial of the creation Heisenberg
operators $\mathfrak a_{-i}(C), \,\, i > 0$ with rational coefficients.

\subsection{\bf Integrality}
\label{subsect_integr}
\par
$\,$

Let $\alpha$ be an arbitrary class in $H^2(X)$.
Fix a smooth irreducible curve $C$ in the surface $X$.
Recall from the previous subsection that the integral operator
$\mathfrak m_{\la, C}$ is a polynomial of the creation Heisenberg
operators $\mathfrak a_{-i}(C), \,\, i > 0$. This enable
us to define $\mathfrak m_{\la, \alpha} \in \End(\mathbb H_X)$ by
replacing the creation operators $\mathfrak a_{-i}(C)$ in
$\mathfrak m_{\la, C}$ by the creation operators $\mathfrak
a_{-i}(\alpha)$ correspondingly, and then define
\begin{eqnarray} \label{def_L_la_a}
[L^\la \alpha] = \mathfrak m_{\la, \alpha} \vac.
\end{eqnarray}
Similarly, we can define the subspaces
$\mathbb H_{n, \alpha}$ and $\mathbb H_\alpha$ of
$\mathbb H_X$ as in (\ref{H_C}).

\begin{lemma} \label{9.14}
Let $\alpha \in H^2(X)$ and $i > 0$. Then, we have
\begin{eqnarray} \label{9.14.1}
\mathfrak a_{-i}(\alpha)[L^\la \alpha] = \sum_{\mu}
a_{\la, \mu} [L^\mu \alpha]
\end{eqnarray}
where the summation is over partitions $\mu$ of $|\la|+i$,
which are obtained as follows:
\begin{enumerate}
\item[{\rm (i)}] add $i$ to a term in $\la$,
say $\la_k$ (possibly $0$), and then

\item[{\rm (ii)}] arrange it in descending order.
\end{enumerate}
The coefficient $a_{\la, \mu}$ is equal to the number of
elements in $\{\ell \,|\, \mu_\ell = \la_k + i\}$.
\end{lemma}
\begin{proof}
Fix a smooth irreducible curve $C$ in $X$.
By the Theorem 9.14 in \cite{Na2},
\begin{eqnarray} \label{9.14.2}
\mathfrak a_{-i}(C)[L^\la C] = \sum_{\mu} a_{\la, \mu} [L^\mu C].
\end{eqnarray}
Define a linear map $\Psi_{C, \alpha}: \mathbb H_C \to
\mathbb H_\alpha$ by sending the basis elements
$\mathfrak a_{-\la}(C) \vac$ of $\mathbb H_C$ to
the elements $\mathfrak a_{-\la}(\alpha) \vac$ in
$\mathbb H_\alpha$. Note that the creation operators
$\mathfrak a_{-j}(C)$ and $\mathfrak a_{-j}(\alpha)$
preserve $\mathbb H_C$ and $\mathbb H_\alpha$ respectively.
Moreover, $\Psi_{C, \alpha} \circ \mathfrak a_{-j}(C) =
\mathfrak a_{-j}(\alpha) \circ \Psi_{C, \alpha}$.
It follows from the definition of $\mathfrak m_{\la, \alpha}$
that $\Psi_{C, \alpha} \circ \mathfrak m_{\la, C} =
\mathfrak m_{\la, \alpha} \circ \Psi_{C, \alpha}$. Thus
\begin{eqnarray} \label{9.14.3}
\Psi_{C, \alpha} \big ([L^\la C] \big )
= \Psi_{C, \alpha} \circ \mathfrak m_{\la, C} \vac
= \mathfrak m_{\la, \alpha} \circ \Psi_{C, \alpha} \vac
= \mathfrak m_{\la, \alpha} \vac = [L^\la \alpha].
\end{eqnarray}
Now applying $\Psi_{C, \alpha}$ to both sides of (\ref{9.14.2}),
we obtain (\ref{9.14.1}).
\end{proof}

\begin{remark}
It is a classical result (cf. \cite{Mac, Na2}) that
\begin{eqnarray} \label{classic}
p_i m_{\lambda} = \sum_{\mu} a_{\la, \mu} m_\mu
\end{eqnarray}
where $a_{\la, \mu}$ is the same as defined in Lemma~\ref{9.14}.
This is compatible with (\ref{Phi_C2}).
\end{remark}

\begin{lemma} \label{-a}
Let $\alpha \in H^2(X)$.
If $[L^\la \alpha]$ is integral for every $\la$,
then so is $[L^\la (-\alpha)]$.
\end{lemma}
\begin{proof}
Fix a smooth irreducible curve $C$ in $X$, and let notations
be the same as in the proof of Lemma~\ref{9.14}.
Recall from (\ref{Phi_C2}) that the linear isomorphism $\Phi_C$
sends $m_\la$ and $p_\la$ to $[L^\la C]$ and
$\mathfrak a_{-\la}(C) \vac$ respectively. Put
\begin{eqnarray} \label{-a.1}
m_\la  = \sum_{|\mu| =|\la|} d_\mu p_\mu
\end{eqnarray}
where $d_\mu \in \mathbb Q$.
Applying $\Psi_{C, \alpha} \circ \Phi_C$ to both sides
and using (\ref{9.14.3}), we obtain
\begin{eqnarray*}
[L^\la \alpha] = \sum_{|\mu| =|\la|}
d_\mu \mathfrak a_{-\mu}(\alpha)\vac.
\end{eqnarray*}
It follows from the definition of $[L^\la(-\alpha)]$ that
\begin{eqnarray} \label{-a.2}
[L^\la(-\alpha)] =\sum_{|\mu| =|\la|} d_\mu (-1)^{\ell(\mu)}
{\mathfrak a}_{-\mu}(\alpha)\vac.
\end{eqnarray}

Next, it is well-known (see \cite{Mac}) that there is
an involution $\omega$ on $\Lambda$ satisfying
\begin{eqnarray*}
\omega (p_\mu) = (-1)^{|\mu| -\ell(\mu)} p_\mu.
\end{eqnarray*}
The so-called forgotten symmetric function associated to
a partition $\mu$ is defined to be $f_\mu =\omega(m_\mu)$.
The symmetric functions $f_\mu, \,\, \mu
\vdash n$ form another $\mathbb Z$-basis of $\Lambda^n$.
In particular, $f_\la$ is an integral linear combination of
the monomial symmetric functions $m_\mu, \,\, \mu \vdash |\la|$.
Thus $\Psi_{C, \alpha} \circ \Phi_C(f_\la)$ is an integral linear
combination of the classes $\Psi_{C, \alpha} \circ \Phi_C(m_\la)
= [L^\mu \alpha], \,\, \mu \vdash |\la|$,
and hence is an integral class. By (\ref{-a.1}),
\begin{eqnarray} \label{-a.3}
\Psi_{C, \alpha} \circ \Phi_C (f_\la)
   &=&\Psi_{C, \alpha} \circ \Phi_C \circ \omega
     \left (\sum_{|\mu| =|\la|} d_\mu p_\mu \right ) \nonumber \\
   &=&(-1)^{|\la|} \cdot \quad \sum_{|\mu| =|\la|} d_\mu
     (-1)^{\ell(\mu)} {\mathfrak a}_{-\mu}(\alpha)\vac.
\end{eqnarray}
Combining (\ref{-a.2}) and (\ref{-a.3}), we conclude that
$[L^\la (-\alpha)]$ is integral as well.
\end{proof}

\begin{lemma} \label{a=a1+a2}
Let $\alpha_1, \alpha_2 \in H^2(X)$,
and $\alpha = \alpha_1 + \alpha_2$. Then,
\begin{eqnarray} \label{a=a1+a2.1}
[L^\la \alpha] = \sum_{(\la^1, \la^2): \la^1 \cup \la^2 = \la}
\mathfrak m_{\la^1, \alpha_1} \mathfrak m_{\la^2, \alpha_2} \vac
\end{eqnarray}
where $(\la^1, \la^2)$ stands for ordered pairs of partitions.
\end{lemma}
\begin{proof}
We start with some notations. For the partitions $\mu$ obtained
from $\la$ as in (\ref{9.14.1}), we shall denote
$\mu = \la\uparrow^{i}$.
If we specify further that such a $\mu$ is obtained from adding
$i$ to a part of $\la$ equal to $j$ (here $j$ is allowed to be
$0$), then we denote $\mu = \la\uparrow_j^{i}$.
Given a partition $\la$, we denote by
$m_k(\la)$ the multiplicity of the parts of $\la$ equal to $k$.
In these notations, the coefficient $a_{\la, \mu}$ in
Lemma~\ref{9.14} for $\mu = \la\uparrow_j^{i}$ is simply
equal to $m_{i+j}(\mu)$. Denote the right-hand-side of
(\ref{a=a1+a2.1}) by $R^\la(\alpha)$.

\medskip\noindent
{\bf Claim.} $\displaystyle{
\mathfrak a_{-i}(\alpha) R^\la(\alpha) = \sum_{\mu
=\la\uparrow^i} a_{\la, \mu} R^\mu (\alpha)}$.
\begin{proof}
The right-hand-side in the Claim is equal to
\begin{eqnarray}  \label{claim.2}
\sum_{\mu =\la\uparrow^i} \,\,
\sum_{(\mu^1, \mu^2): \mu^1 \cup \mu^2 = \mu}
    a_{\la, \mu} \mathfrak m_{\mu^1, \alpha_1} \mathfrak
m_{\mu^2, \alpha_2}
    \vac.
\end{eqnarray}

On the other hand, by Lemma~\ref{9.14},
the left-hand-side in the claim is
\begin{eqnarray}
  && \mathfrak a_{-i}(\alpha) R^\la(\alpha)  \nonumber \\
  &=& \sum_{(\la^1, \la^2): \la^1 \cup \la^2 = \la}
    ( \mathfrak m_{\la^2, \alpha_2}
    \mathfrak a_{-i}(\alpha_1) \mathfrak m_{\la^1,
\alpha_1} \vac
     +
    \mathfrak m_{\la^1, \alpha_1}
    \mathfrak a_{-i}(\alpha_2) \mathfrak m_{\la^2,
\alpha_2} \vac )
    \nonumber  \\
  &=&  \sum_{(\la^1, \la^2): \la^1 \cup \la^2 = \la} \;
      \sum_{\rho^1 =\la^1 \uparrow^i}  a_{\la^1, \rho^1}
\,\, \mathfrak m_{\rho^1, \alpha_1}
    \mathfrak m_{\la^2, \alpha_2}\vac   \label{sum1}\\
    && + \sum_{(\la^1, \la^2): \la^1 \cup \la^2 = \la} \;
    \sum_{\rho^2 =\la^2 \uparrow^i}
    a_{\la^2, \rho^2} \,\, \mathfrak m_{\la^1, \alpha_1}
      \mathfrak m_{\rho^2, \alpha_2}  \vac.  \label{sum2}
\end{eqnarray}
Note that the partitions $\rho^1 \cup \la^2$ and
$\la^1 \cup \rho^2$ associated to $\la^1, \la^2, \rho^1$
and $\rho^2$ appearing above are of
the form $\la \uparrow^i$. Thus the same type of terms appearing
on both sides of the Claim. It remains to identify the
coefficients of a given term.

Fix $\mu =\la \uparrow^i_j$ for some part $j$ of $\la$, and fix
$\mu^1,\mu^2$ such that $\mu^1 \cup \mu^2 =\mu$. From the above
computation, the contributions to the term
$\mathfrak m_{\mu^1, \alpha_1} \mathfrak m_{\mu^2, \alpha_2} \vac$
in the left-hand-side of the Claim come from two places:
the term in (\ref{sum1}) for $\rho^1 =\mu^1, \la^2 =\mu^2$
whose coefficient is $m_{i+j}(\mu^1)$, and the term in (\ref{sum2})
for $\la^1 =\mu^1, \rho^2 =\mu^2$ whose coefficient is
$m_{i+j}(\mu^2)$. Therefore in view of (\ref{claim.2}),
the coefficients of the term $\mathfrak m_{\mu^1, \alpha_1}
\mathfrak m_{\mu^2, \alpha_2}\vac$
in both sides of the Claim coincide thanks to
\begin{eqnarray}   \label{extra}
m_{i+j} (\mu^1) +m_{i+j} (\mu^2) = m_{i+j} (\mu^1 \cup \mu^2)
= m_{i+j} (\mu) =a_{\la,\mu}.
\end{eqnarray}
This completes the proof of the above Claim.
\end{proof}

Next, we continue the proof of (\ref{a=a1+a2.1}) by using
induction on $n$ and the reverse dominance ordering of
partitions $\la$ of $n$ (we refer to p.7 of \cite{Mac}
for the definition of the dominance ordering).
For $n=1$, formula (\ref{a=a1+a2.1}) is clear.
Assume that formula (\ref{a=a1+a2.1}) holds for
all partitions of size less than $n$.

For $\la =(n)$, (\ref{a=a1+a2.1}) holds since $\mathfrak m_{\la,
\alpha} = \mathfrak a_{-n}(\alpha)$ and $[L^\la \alpha] =
\mathfrak a_{-n}(\alpha) \vac$.

For a general partition $\la$ of $n$ with a part equal to,
say $i$, we denote by $\widetilde{\la}$ the partition obtained
from $\la$ with a part equal to $i$ removed.
Now replacing $\la$ in (\ref{9.14.1}) and the above Claim by
$\widetilde{\la}$ respectively, we obtain
\begin{eqnarray}  \label{terms}
\mathfrak a_{-i}(\alpha)[L^{\widetilde{\la}} \alpha]
  &=& \sum_{\mu=\widetilde{\la}\uparrow^i}
a_{\widetilde{\la}, \mu} [L^\mu \alpha],  \label{term1} \\
\mathfrak a_{-i}(\alpha) R^{\widetilde{\la}}(\alpha)
  &=& \sum_{\mu=\widetilde{\la}\uparrow^i}
a_{\widetilde{\la}, \mu} R^\mu (\alpha). \label{term2}
\end{eqnarray}
Note that $\la$ appears among the above $\mu$'s as the maximum
in the reverse dominance ordering. By induction hypothesis,
the left-hand-sides of (\ref{term1}) and (\ref{term2}) coincide,
and all the terms on the right-hand-sides of (\ref{term1}) and
(\ref{term2}) involving $\mu$ not equal to $\la$ coincide.
Thus (\ref{a=a1+a2.1}) follows since $a_{\widetilde{\la}, \la}
= m_i(\la) \neq 0$.
\end{proof}

\begin{theorem} \label{div}
For every divisor $\alpha$ on $X$ and every partition $\la$,
$[L^\la \alpha]$ is an integral class and $\mathfrak m_{\la,
\alpha} \in \End(\mathbb H_X)$ is an integral operator.
\end{theorem}
\begin{proof}
Every divisor $\alpha$ can be written as $C_1 - C_2$ for some very
ample divisors $C_1$ and $C_2$. Represent $C_1$ and $C_2$ by
smooth irreducible curves. Now we see from the discussions in
Subsection~\ref{subsect_c}, Lemma~\ref{-a} and Lemma~\ref{a=a1+a2}
that $[L^\la \alpha]$ is an integral class. By
Proposition~\ref{prop_int}, $\mathfrak m_{\la, \alpha} \in
\End(\mathbb H_X)$ is an integral operator.
\end{proof}

\subsection{\bf Unimodularity}
\label{subsect_unimodu}
\par
$\,$

Our goal in this subsection is to prove the following.

\begin{theorem} \label{unimod}
Let $\alpha_1, \cdots, \alpha_k \in H^2(X)$ be linearly
independent classes, and let $M_{\underline{\alpha}}$ be the
intersection matrix of $\alpha_1, \cdots, \alpha_k$.
Fix a positive integer $n$.
Let $M_{n, \underline{\alpha}}$ be the intersection matrix
of the classes in $H^{2n}(\Xn)$:
\begin{eqnarray}  \label{unimod.1}
\mathfrak m_{\la^1, \alpha_1} \cdots \mathfrak
m_{\la^k, \alpha_k} \vac, \quad |\la^1| + \ldots + |\la^k| = n.
\end{eqnarray}
If $\det M_{\underline{\alpha}} = \pm 1$, then we have
$\det M_{n, \underline{\alpha}} = \pm 1$ as well.
\end{theorem}

First of all, we begin with some linear algebra preparations.
Let $V$ be a $k$-dimensional complex vector space
with a symmetric bilinear form
\begin{eqnarray*}
(-, -): V\times V \to \C.
\end{eqnarray*}
Fix a linear basis ${\bf v} =\{v_1, \ldots, v_k\}$ of $V$, and let
$M_{\bf v} = \big ( (v_i, v_j) \big )_{1 \le i,j \le k}$
be the matrix formed by the pairings $(v_i, v_j)$.
The following is elementary.

\begin{lemma} \label{lin_alg_lma1}
Let $\tilde{\bf v} =\{\tilde{v}_1, \ldots, \tilde{v}_k\}$ be
another basis of $V$. Let $T$ be the transition matrix from
$\tilde{\bf v}^t$ to ${\bf v}^t$ (i.e., ${\bf v}^t =
T \tilde{\bf v}^t$). Then,
$M_{\bf v} = T M_{\tilde{\bf v}} T^t$. \hfill \qed
\end{lemma}

For a given $n \ge 0$, the symmetric power $S^n(V)$
admits a monomial basis
\begin{eqnarray}  \label{mono_basis}
\{v_{i_1} \cdots v_{i_n}|\,\,
 1 \le i_1 \le \ldots \le i_n \le k \}.
\end{eqnarray}
Define a bilinear form, still denoted by $(-, -)$,
on $S^n(V)$ by letting
\begin{eqnarray*}
( v_{i_1} \cdots v_{i_n}, v_{j_1} \cdots v_{j_n}) =
\sum_{\sigma \in S_n} \prod_{a=1}^n  ( v_{i_a}, v_{ j_{\sigma (a)}})
\end{eqnarray*}
where $S_n$ denotes the $n$-th symmetric group.
Denote by $M_{n, \bf v}$ the matrix of the pairings of the
monomial basis elements of $S^n(V)$
(so we have $M_{1, \bf v} = M_{\bf v}$).

\begin{lemma} \label{lin_alg_lma2}
For some constant $c(n,k)$, we have
\begin{eqnarray} \label{poly}
{\det M_{n, \bf v} = c(n,k) \cdot \left ( \det M_{\bf v} \right
)^{{n+k-1 \choose k}}}.
\end{eqnarray}
\end{lemma}
\par\noindent
\begin{proof}
First we prove the formula under the assumption that
\begin{eqnarray} \label{assume}
\det M_{{\bf v}_i} \ne 0 \quad \text{ for }1 \le i \le k-1
\end{eqnarray}
where ${\bf v}_i := \{v_1, \ldots, v_i \}$ for $1 \le i \le k$
(so ${\bf v}_k = {\bf v}$).

Use induction on $k$. When $k=1$, our formula is clearly true.

Fix $k > 1$. Since $\det M_{{\bf v}_{k-1}} \ne 0$, there is a
vector $\tilde{v}_k = e_1v_1 + \ldots + e_{k-1}v_{k-1} + v_k$
such that $e_1, \ldots, e_{k-1} \in \C$ and $(v_1, \tilde{v}_k) =
\ldots = (v_{k-1}, \tilde{v}_k) = 0$. Let
\begin{eqnarray*}
\tilde{\bf v} =\{v_1, \ldots, v_{k-1}, \tilde{v}_k \}.
\end{eqnarray*}
Then the transition matrix from the basis $\tilde{\bf v}^t$ to the
basis ${\bf v}^t$ is lower triangular with all diagonal entries
being $1$. By Lemma~\ref{lin_alg_lma1}, we have
\begin{eqnarray}  \label{lin_alg_lma2.1}
\det M_{\bf v}  = \det M_{\tilde{\bf v}} = \det M_{{\bf v}_{k-1}}
\cdot (\tilde{v}_k, \tilde{v}_k).
\end{eqnarray}

Next, we apply Lemma~\ref{lin_alg_lma1} to $S^n(V)$ for the two
monomial bases:
\begin{eqnarray*}
\{ v_1^{i_1} \cdots v_{k-1}^{i_{k-1}} v_k^{i_k}
   |\,\, i_1 + \ldots + i_{k-1} + i_k =n \},  \\
\{ v_1^{i_1} \cdots v_{k-1}^{i_{k-1}}
   \tilde{v}_k^{i_k} |\,\, i_1 + \ldots + i_{k-1} + i_k =n \}.
\end{eqnarray*}
Since the transition matrix from the second basis to the first
basis is lower triangular with all diagonal entries being $1$, we
conclude that
\begin{eqnarray}  \label{lin_alg_lma2.2}
\det M_{n, \bf v} = \det M_{n, \tilde{\bf v}}.
\end{eqnarray}

Now note that $(v_1^{i_1} \cdots v_{k-1}^{i_{k-1}}
\tilde{v}_k^{i_k}, \,\, v_1^{j_1} \cdots v_{k-1}^{j_{k-1}}
\tilde{v}_k^{j_k})$ is equal to
$$(v_1^{i_1} \cdots v_{k-1}^{i_{k-1}}, \,\, v_1^{j_1} \cdots
v_{k-1}^{j_{k-1}}) \cdot \delta_{i_k, j_k} \cdot
i_k! \cdot (\tilde{v}_k, \tilde{v}_k)^{i_k}.$$
It follows that $M_{n, \tilde{\bf v}} = \text{diag}
\big ( \cdots, \,\, M_{n-m, {\bf v}_{k-1}} \cdot
m! \cdot (\tilde{v}_k, \tilde{v}_k)^m, \,\, \cdots \big )$
where $m$ runs from $0$ to $n$.
The matrix $M_{n-m, {\bf v}_{k-1}}$ has
${n-m+k-2 \choose k-2}$ rows. So
\begin{eqnarray*}
\det M_{n, \tilde{\bf v}} = c_1(n,k) \cdot \prod_{m=0}^n \left (
\det M_{n-m, {\bf v}_{k-1}} \cdot (\tilde{v}_k, \tilde{v}_k)^{m
\cdot {n-m+k-2 \choose k-2}} \right )
\end{eqnarray*}
for some constant $c_1(n,k)$. Applying induction to
$\det M_{n-m, {\bf v}_{k-1}}$ yields
\begin{eqnarray}  \label{lin_alg_lma2.3}
   \det M_{n, \tilde{\bf v}}
&=&c(n,k) \cdot \prod_{m=0}^n \left (
   \big (\det M_{{\bf v}_{k-1}} \big )^{n-m+k-2 \choose k-1}
   \cdot (\tilde{v}_k, \tilde{v}_k)^{m \cdot {n-m+k-2 \choose k-2}}
  \right )    \nonumber \\
&=&c(n,k) \cdot
   \big (\det M_{{\bf v}_{k-1}} \big )^{\sum_{m=0}^n
   {n-m+k-2 \choose k-1}} \cdot
  (\tilde{v}_k, \tilde{v}_k)^{\sum_{m=0}^n
  {m \cdot {n-m+k-2 \choose k-2}}}    \nonumber \\
&=&c(n,k) \cdot \big (\det M_{{\bf v}_{k-1}} \cdot
  (\tilde{v}_k, \tilde{v}_k)\big )^{n+k-1 \choose k}
\end{eqnarray}
for some constant $c(n,k)$, where we have used the combinatorial
identities:
\begin{eqnarray}  \label{comb_id}
\sum_{m=0}^n {n-m+k-2 \choose k-1} = \sum_{m=0}^n {m \cdot
{n-m+k-2 \choose k-2}} = {n+k-1 \choose k}
\end{eqnarray}
Indeed, all the three terms in (\ref{comb_id}) compute the
dimension of the space of degree-$(n-1)$ homogeneous polynomials
in $(k+1)$-variables.

Combining (\ref{lin_alg_lma2.2}), (\ref{lin_alg_lma2.3})
and (\ref{lin_alg_lma2.1}), we conclude that
\begin{eqnarray*}
  \det M_{n, \bf v}
= c(n,k) \cdot \big (\det M_{{\bf v}_{k-1}} \cdot
  (\tilde{v}_k, \tilde{v}_k)\big )^{n+k-1 \choose k}
= c(n,k) \cdot \big ( \det M_{\bf v} \big )^{n+k-1 \choose k}.
\end{eqnarray*}

Finally, we come to the general case. Assume that
the bilinear form $(-, -)$ on $V$ is not identically zero
(otherwise the lemma trivially holds).
Set $z_{ij} =(v_i, v_j)$ for $1 \le i \le j \le k$.
Then both sides of (\ref{poly}) are easily seen to be
polynomials in the variables $z_{ij}, \,\, 1 \le i \le j \le k$.
The proof above under the assumption
(\ref{assume}) implies that the (polynomial) identity
(\ref{poly}) holds for a Zariski open subset
$(z_{ij})_{1 \le i \le j \le k}$ of $\C^{k(k+1)/2}$, and
thus it holds for an arbitrary $(z_{ij})_{1 \le i \le j \le k}$.
\end{proof}

For $r \ge 1$, we denote by $V[r]$ a copy of $V$ with
bilinear form
\begin{eqnarray} \label{()r}
(-,-)_r = (-1)^{r-1}r \cdot (-,-).
\end{eqnarray}
We shall denote by $(-,-)_r$ as well the induced bilinear form
on the symmetric power $S^*(V[r])$.
In particular, $V[1] =V$ with $(-,-)_1 =  (-,-).$

For a partition $\mu =(r^{m_r})_{r \ge 1}$,
we form the vector space
$S^\mu V := \bigotimes_{r \ge 1} S^{m_r} (V[r])$
with a bilinear form $(-,-)_\mu$ given by
$(\otimes_r u_{I_r}, \otimes_r v_{J_r})_\mu
 =\prod_r (u_{I_r},   v_{J_r})_r$
for $u_{I_r}, v_{J_r} \in S^{m_r} (V[r])$.
Denote by $M_{\mu, \bf v}$ the matrix of the pairings of
the induced monomial basis for $S^\mu V$
(see (\ref{mono_basis}) for the monomial basis of $S^n(V)$).

\begin{lemma} \label{iden}
For some constant $c(\mu,k)$ and some integer $d(\mu,k) \ge 1$, we
have
\begin{eqnarray*}
\displaystyle{\det M_{\mu, \bf v} = c(\mu,k) \cdot \big ( \det
M_{\bf v} \big )^{d(\mu,k)}}.
\end{eqnarray*}
\end{lemma}
\begin{proof}
Follows immediately from our definitions and
Lemma~\ref{lin_alg_lma2}.
\end{proof}

Now let $\alpha_1, \ldots, \alpha_k \in H^2(X)$ be from
Theorem~\ref{unimod}.
Let $V = \C \alpha_1 \oplus \cdots \oplus \C \alpha_k$ with
the pairing induced from the one on $H^*(X)$.
Let $\mathbb H_{n, V} \subset H^*(\Xn)$ be the $\C$-linear
span of the classes (\ref{unimod.1}). Then,
$\mathbb H_{n, V}$ has two linear bases:
\begin{eqnarray}
\mathfrak m_{n, \underline{\alpha}}
  &:=& \{\mathfrak m_{\la^1, \alpha_1}
  \cdots \mathfrak m_{\la^k, \alpha_k} \vac \},
  \label{m_n_a}  \\
\mathfrak a_{n, \underline{\alpha}}
  &:=& \{\mathfrak a_{-\la^1}(\alpha_1)
  \cdots \mathfrak a_{-\la^k}(\alpha_k) \vac \} \label{a_n_a}
\end{eqnarray}
where $\la^1, \ldots, \la^k$ are partitions satisfying
$|\la^1| + \cdots + |\la^k| = n$.

Take an orthogonal basis $\{\beta_1, \cdots, \beta_k \}$ of
$V$ with $(\beta_i, \beta_i) = \pm 1$ for every $i$.
We define similarly bases
$\mathfrak m_{n, \underline{\beta}}$ and
$\mathfrak a_{n, \underline{\beta}}$ for $\mathbb H_{n, V}$.
Observe that the transition matrix $B$ from
the basis $\mathfrak m_{n, \underline{\alpha}}$ to
$\mathfrak a_{n, \underline{\alpha}}$ is independent of
the basis $\{\alpha_1, \cdots, \alpha_k \}$ of $H^2(X)$,
that is, it is the same as the transition matrix from
the basis $\mathfrak m_{n, \underline{\beta}}$ to
$\mathfrak a_{n, \underline{\beta}}$.

Denote by $A$ the transition matrix from the basis
$\mathfrak a_{n, \underline{\beta}}$ to
$\mathfrak a_{n, \underline{\alpha}}$.
Denote by $M_{n, \underline{\beta}}$ the intersection matrix
of the basis $\mathfrak m_{n, \underline{\beta}}$.
Since the transition matrix from
$\mathfrak m_{n, \underline{\beta}}$ to
$\mathfrak m_{n, \underline{\alpha}}$ is $B^{-1} A B$,
we obtain from Lemma~\ref{lin_alg_lma1} the following.

\begin{lemma} \label{linear}
$M_{n, \underline{\alpha}} = (B^{-1} A B) \,\,
M_{n, \underline{\beta}} \,\,(B^{-1} A B)^t$. \qed
\end{lemma}

Next, we prove a special case of Theorem~\ref{unimod} for
the basis $\{\beta_1, \cdots, \beta_k \}$.

\begin{lemma} \label{det=one}
$\det M_{n, \underline{\beta}} = \pm 1$.
\end{lemma}
\begin{proof}
As complex vector spaces, we have a linear isomorphism:
\begin{eqnarray} \label{det=one.0}
\mathbb H_{n,V} \cong \bigoplus_{\underline{n}} \bigotimes_{i=1}^k
\mathbb H_{n_i, \C \beta_i}
\end{eqnarray}
where the sum is over $\underline{n} = (n_1, \cdots, n_k) \in
(\Z_{\ge 0})^k$ such that $n_1 + \ldots + n_k =n$.
Denote by $M_{n_i,\beta_i}$ the matrix of pairings among
the basis $\big \{\mathfrak m_{\la^i, \beta_i} \vac \, | \,\, \la^i
\vdash n_i \big \}$ for $\mathbb H_{n_i, \C \beta_i}$. For
given $1 \le i \le k$ and $\lambda^i$, the operator $\mathfrak
m_{\la^i, \beta_i}$ is a linear combination of the operators
$\mathfrak a_{- \mu^i}(\beta_i)$, $\mu^i \vdash |\la^i|.$ Since
the classes $\beta_1, \cdots, \beta_k$ are orthogonal,
we see from (\ref{eq:heis}) that
$\big ( \mathfrak m_{\la^1, \beta_1} \cdots
   \mathfrak m_{\la^k, \beta_k} \vac,
   \mathfrak m_{\mu^1, \beta_1} \cdots
   \mathfrak m_{\mu^k, \beta_k} \vac\big )$ is equal to
\begin{eqnarray*}
\prod_{i=1}^k \left ( \delta_{|\la^i|, |\mu^i|} \,\,
   \cdot \,\, \big ( \mathfrak m_{\la^i, \beta_i}\vac,
   \mathfrak m_{\mu^i, \beta_i} \vac \big ) \right ).
\end{eqnarray*}
This orthogonality together with (\ref{det=one.0}) implies
by standard linear algebra that
$$\det M_{n, \underline{\beta}} =\prod_{\underline{n}} \prod_{i=1}^k
(\det M_{n_i,\beta_i})^{\prod_{1 \le j \le k, j \ne i}
{\rm dim} \, \mathbb H_{n_j, \C \beta_j}}.
$$
Thus it suffices to prove our lemma for $k = 1$.

Put $\beta = \beta_1$ and $V =\C\beta$. By (\ref{m_n_a}) and (\ref{a_n_a}),
$\mathbb H_{n, V}$ has two linear bases:
\begin{eqnarray*}
\mathfrak m_{n, \underline{\beta}}
  = \big \{\mathfrak m_{\la, \beta}\vac|\,\, \la \vdash n \big \},\qquad
\mathfrak a_{n, \underline{\beta}}
  = \big \{\mathfrak a_{-\la}(\beta)\vac|\,\, \la \vdash n  \big \}.
\end{eqnarray*}
Since $(\beta,\beta)=\pm 1$, we obtain by (\ref{a_n_pm}) and
(\ref{eq:heis}) that
\begin{eqnarray} \label{det=one.01}
\big (\mathfrak a_{-\la}(\beta)\vac, \mathfrak a_{-\mu}(\beta)\vac
\big ) = \pm \delta_{\la, \mu} \cdot {\mathfrak z}_\la \cdot
(\beta, \beta)^{\ell(\la)} = \pm \delta_{\la, \mu} \cdot
{\mathfrak z}_\la.
\end{eqnarray}
Recall that the transition matrix from the basis $\mathfrak a_{n,
\underline{\beta}}$ to $\mathfrak m_{n, \underline{\beta}}$ is
$B^{-1}$. So we conclude from Lemma~\ref{lin_alg_lma1}
and (\ref{det=one.01}) that
$M_{n, \underline{\beta}} \, = \, B^{-1} \cdot \text{diag}(\cdots,
\pm {\mathfrak z}_\la, \cdots) \cdot (B^{-1})^t$
where $\la$ runs over all partitions of $n$. Hence we have
\begin{eqnarray} \label{det=one.1}
\det M_{n, \underline{\beta}} \, = \, \pm (\det B)^{-2} \cdot
\prod_{\la \vdash n} {\mathfrak z}_\la
\end{eqnarray}

By (\ref{Phi_C2}), Lemma~\ref{9.14}, and (\ref{classic}),
we see that the transition matrix from the basis $\{p_\la| \,\,
\la \vdash n \}$ of $\Lambda_{\mathbb Q}^n \subset \Lambda_{\mathbb Q}$ to
the basis $\{m_\la| \,\, \la \vdash n\}$ is $B^{-1}$ as well.
Introduce the standard bilinear form $(-, -)$ on
$\Lambda_{\mathbb Q}^n$ (cf. \cite{Mac}) by letting
\begin{eqnarray} \label{pairing_Lam}
(p_\la, p_\mu) = \delta_{\la, \mu} \cdot {\mathfrak z}_\la.
\end{eqnarray}
The matrix $M_n$ formed by the pairings of the monomial symmetric
functions $m_\la, \,\, \la \vdash n$ is unimodular, and thus $\det
M_n = \pm 1$. By Lemma~\ref{lin_alg_lma1} and (\ref{pairing_Lam}),
\begin{eqnarray*}
M_n \, = \, B^{-1} \cdot  \text{diag}(\cdots, {\mathfrak z}_\la,
\cdots)  \cdot (B^{-1})^t.
\end{eqnarray*}
where $\la$ runs over all partitions of $n$. It follows immediately that
\begin{eqnarray} \label{det_Tn}
(\det B)^{-2} \cdot \prod_{\la \vdash n} {\mathfrak z}_\la = \det
M_n = \pm 1.
\end{eqnarray}
Combining this with (\ref{det=one.1}), we finally obtain
$\det M_{n, \underline{\beta}} = \pm 1$.
\end{proof}

For $\mu \vdash n$, let $\mathbb H_{\mu, \underline{\alpha}} \subset
\mathbb H_{n, V}$ be the span of
$\mathfrak a_{-\mu^1}(\alpha_1)
\cdots \mathfrak a_{-\mu^k}(\alpha_k) \vac$
where $\mu^1, \ldots, \mu^k$ are partitions such that
$\mu^1 \cup \cdots \cup \mu^k =\mu$. Note that
$\mathbb H_{\mu, \underline{\alpha}} \subset
\mathbb H_{n, V} \subset H^{2n}(\Xn)$.
So $\mathbb H_{\mu, \underline{\alpha}}$ carries a pairing
induced from the one on $H^{2n}(\Xn)$.

\begin{lemma} \label{isometry}
For every partition $\mu$,
we have an isometry
$\mathbb H_{\mu, \underline{\alpha}} \cong S^{\mu} (V)$.
\end{lemma}
\begin{proof}
Follows from the definitions, (\ref{()r}) and the commutation
relation (\ref{eq:heis}).
\end{proof}

\begin{lemma}  \label{simple}
$(\det A)^2 = \pm 1$.
\end{lemma}
\begin{proof}
By the Heisenberg algebra commutation relation (\ref{eq:heis}),
$\mathfrak a_{-\mu^1}(*) \cdots \mathfrak a_{-\mu^k}(*)$
is orthogonal to $\mathfrak a_{-\nu^1}(*) \cdots
\mathfrak a_{-\nu^\ell}(*) \vac$ unless
$\mu^1 \cup \cdots \cup \mu^k
= \nu^1 \cup \cdots \cup \nu^\ell$.
Here $*$'s denote unspecified classes in $H^2(X)$. Thus,
we have an orthogonal direct sum
$$\mathbb H_{n, V} =
\bigoplus_{ \mu \vdash n} \mathbb H_{\mu, \underline{\alpha}}.$$
By Lemma~\ref{isometry}, the intersection matrix
$\widetilde{M}_{n, \underline{\alpha}}$ of the basis
$\mathfrak a_{n, \underline{\alpha}}$ for $\mathbb H_{n, V}$
is given by the diagonal block matrix whose diagonal
consists of $M_{\mu, \underline{\alpha}}$, $\mu \vdash n$. So
\begin{eqnarray} \label{simple.1}
\det \widetilde{M}_{n, \underline{\alpha}}
= \prod_{\mu \vdash n} \det M_{\mu, \underline{\alpha}}.
\end{eqnarray}

Similarly, by repeating the above with the $\alpha$'s
replaced by the $\beta$'s, we obtain
\begin{eqnarray} \label{simple.2}
\det \widetilde{M}_{n, \underline{\beta}} = \prod_{\mu \vdash n}
\det M_{\mu, \underline{\beta}}.
\end{eqnarray}

By assumptions, $\det M_{\underline{\alpha}} =\pm 1$ and
$\det M_{\underline{\beta}} = \pm 1$. Thus by Lemma~\ref{iden},
we have $\det M_{\mu, \underline{\alpha}}
= \pm \det M_{\mu, \underline{\beta}}$ for every $\mu$.
Therefore, we see from (\ref{simple.1}) and
(\ref{simple.2}) that
\begin{eqnarray} \label{simple.3}
\det \widetilde{M}_{n, \underline{\alpha}}
= \pm \det \widetilde{M}_{n, \underline{\beta}}.
\end{eqnarray}
By Lemma~\ref{lin_alg_lma1} and the definition of the matrix $A$,
we have $\widetilde{M}_{n, \underline{\alpha}}
= A \, \widetilde{M}_{n, \underline{\beta}} A^t$.
Combining this with (\ref{simple.3}) yields $(\det A)^2 = \pm 1$.
\end{proof}

Finally, Theorem~\ref{unimod} follows from Lemmas \ref{linear},
\ref{det=one}, and \ref{simple}.

\section{\bf Integral bases for the cohomology of
Hilbert schemes}
\label{sect_app}

\subsection{\bf Structure of integral bases for $H^*(\Xn)$}
\label{subsect_str}
\par
$\,$

\begin{definition} \label{def_H'}
For $n \ge 0$, define ${\mathbb H}_{n, X}'$ to be the linear
subspace of $\mathbb H_{n, X} := H^*(\Xn)$ spanned by all
Heisenberg monomial classes involving only the creation operators
$\mathfrak a_{-j}(\alpha)$ where $j > 0$ and $\alpha \in H^1(X)
\oplus H^2(X) \oplus H^3(X)$.
\end{definition}

Parallel to (\ref{Ai}), every class $A \in H^*(\Xn)$ can
be written as
\begin{eqnarray} \label{A_la_mu}
A = \sum_{\la, \mu} \frac{1}{{\mathfrak z}_\la}
\cdot \mathfrak a_{-\la}(1)\mathfrak a_{-\mu}(x)(A_{\la, \mu})
\end{eqnarray}
where $A_{\la, \mu} \in {\mathbb H}_{n-|\la|-|\mu|, X}'$.
The following is an analogue of Lemma~\ref{int_Ai}.

\begin{lemma} \label{int_A_la_mu}
Let $A \in H^*(\Xn)$ be expressed as in (\ref{A_la_mu}). Then,
\begin{enumerate}
\item[{\rm (i)}] the classes $A_{\la, \mu}$ are uniquely
determined by $A$;

\item[{\rm (ii)}] $A$ is integral if and only if all
the classes $A_{\la, \mu}$ are integral.
\end{enumerate}
\end{lemma}
\begin{proof}
Our argument is similar to the one used in the proof of
Lemma~\ref{int_Ai}. First of all, the operator $1/\mathfrak z_\la
\cdot \mathfrak a_{-\la}(1)$ are integral by Lemma~\ref{int_a_la}.
Since $\mathfrak a_{-\mu}(x)$ is also integral, we see that the
class $A$ is integral if all the classes $A_{\la, \mu}$ are
integral. So it remains to prove (i) and the ``only if" part of
(ii).

Next, let $n_0$ be the maximal integer such that $A_{\la, \mu} \ne
0$ for some partitions $\la, \mu$ with $|\la|+|\mu|= n_0$. If
$\la$ and $\mu$ are partitions with $|\la|+|\mu|= n_0$, then
applying the adjoint operator $\big ( 1/{\mathfrak z}_\mu \cdot
\mathfrak a_{-\mu}(1)\mathfrak a_{-\la}(x) \big )^\dagger$ to both
sides of (\ref{A_la_mu}) yields:
\begin{eqnarray*}
\big ( 1/{\mathfrak z}_\mu \cdot \mathfrak a_{-\mu}(1)
   \mathfrak a_{-\la}(x) \big )^\dagger(A)
= (-1)^{|\la|-\ell(\la) + |\mu|-\ell(\mu)} \cdot A_{\la, \mu}.
\end{eqnarray*}
So the class $A_{\la, \mu}$ is uniquely determined by $A$.
Moreover, by Lemma~\ref{int_Heis}~(i), if $A$ is
an integral class, then so is the class $A_{\la, \mu}$.

Finally, repeating the above process to the class
\begin{eqnarray*}
    A'
&:=&A - \sum_{|\la|+|\mu|= n_0} \frac{1}{{\mathfrak z}_\la}
    \cdot \mathfrak a_{-\la}(1)
    \mathfrak a_{-\mu}(x)(A_{\la, \mu})  \\
&= &\sum_{|\la|+|\mu| < n_0} \frac{1}{{\mathfrak z}_\la}
    \cdot \mathfrak a_{-\la}(1)
    \mathfrak a_{-\mu}(x)(A_{\la, \mu}),
\end{eqnarray*}
we conclude that (i) and the ``only if" part of (ii) hold
for all the classes $A_{\la, \mu}$.
\end{proof}

\begin{proposition} \label{str}
Assume that for each fixed $k \ge 0$, the cohomology classes
$B_{k, i}, i \in \mathcal I_k$ form an integral basis of
${\mathbb H}_{k, X}'$, where $\mathcal I_k$ is an index set
depending on $k$.
Then, an integral basis of $H^*(\Xn)$ consists of the classes:
\begin{eqnarray} \label{str.1}
\frac{1}{{\mathfrak z}_\la}
\cdot \mathfrak a_{-\la}(1)\mathfrak a_{-\mu}(x)B_{k, i}, \qquad
|\la|+|\mu| + k = n, \,\, i \in \mathcal I_k.
\end{eqnarray}
\end{proposition}
\begin{proof}
Since $B_{k, i}$ are integral classes and the operators
$1/{\mathfrak z}_\la \cdot \mathfrak a_{-\la}(1)
\mathfrak a_{-\mu}(x)$ are integral,
the cohomology classes (\ref{str.1}) are integral.
By Lemma~\ref{int_A_la_mu}~(ii), every integral class in
$H^*(\Xn)$ is an integral linear combination of
the classes (\ref{str.1}).
Since the dimension of $H^*(\Xn)$ is equal to the number of
elements in (\ref{str.1}), we conclude that an integral basis
of $H^*(\Xn)$ is given by the classes (\ref{str.1}).
\end{proof}

This proposition enables us to write down an integral basis of
$H^*(\Xn)$ whenever an integral basis of the subspace
${\mathbb H}_{k, X}'$ is known for every $k \le n$.

\subsection{\bf Surfaces $X$ with
$H^1(X; \mathcal O_X) = H^2(X; \mathcal O_X) = 0$}
\label{subsect_pg_q}
\par
$\,$

\begin{theorem} \label{thm_pg_q}

Let $X$ be a projective surface with $H^1(X; \mathcal O_X) =
H^2(X; \mathcal O_X) = 0$. Let $\alpha_1, \cdots, \alpha_k$ be an
integral basis of $H^2(X)$. Then the following classes
\begin{eqnarray} \label{thm_pg_q.1}
\frac{1}{{\mathfrak z}_\la} \cdot \mathfrak a_{-\la}(1) \mathfrak
a_{-\mu}(x) \mathfrak m_{\nu^1, \alpha_1} \cdots \mathfrak
m_{\nu^k, \alpha_k} \vac, \quad |\la|+|\mu|+ \sum_{i=1}^k |\nu^i|
= n
\end{eqnarray}
are integral, and furthermore, they form an integral basis for
$H^*(\Xn; \mathbb Z)/\text{\rm Tor}$.
\end{theorem}
\begin{proof}
By Proposition~\ref{str}, it remains to show that for each
fixed $\ell \ge 0$, an integral basis of ${\mathbb H}_{\ell, X}'$
is given by the integral cohomology classes:
\begin{eqnarray} \label{thm_pg_q.2}
\mathfrak m_{\nu^1, \alpha_1}
\cdots \mathfrak m_{\nu^k, \alpha_k} \vac,
\quad \sum_{i=1}^k |\nu^i| = \ell.
\end{eqnarray}

Since $H^1(X; \mathcal O_X) = H^2(X; \mathcal O_X) = 0$,
we have $H^2(X; \Z) \cong \text{\rm Pic}(X)$ and
$H^1(X) = H^3(X) = 0$. Hence the integral classes
$\alpha_1, \cdots, \alpha_k$ are divisors.
By Theorem~\ref{div}, all the classes in (\ref{thm_pg_q.2})
are integral. By our assumption, the intersection matrix of
the classes $\alpha_1, \cdots, \alpha_k$ is unimodular.
It follows from Theorem~\ref{unimod} that the intersection
matrix of the classes (\ref{thm_pg_q.2}) is unimodular as well.
Since $H^1(X) = H^3(X) = 0$, the dimension of the space
${\mathbb H}_{\ell, X}'$ is precisely equal to the number of
classes in (\ref{thm_pg_q.2}).
Therefore, the classes (\ref{thm_pg_q.2})
form an integral basis of ${\mathbb H}_{\ell, X}'$.
\end{proof}

We conjecture that the cohomology class $[L^\la \alpha]$ is
integral whenever $\alpha \in H^2(X)$ is an integral class.
If this conjecture is true, then the statement in
Theorem~\ref{thm_pg_q} will be valid for
every projective surface $X$ with vanishing odd cohomology.

\begin{remark} \label{rmk_geom}
Let $C_1, \cdots, C_k$ be smooth irreducible curves in $X$ such
that any two of them intersect transversely and no three of them
intersect. Let $\nu^1, \cdots, \nu^k$ be partitions with $|\nu^1|+
\cdots + |\nu^k| = \ell$. Then, $\mathfrak m_{\nu^1, C_1} \cdots
\mathfrak m_{\nu^k, C_k} \vac \in H^{2 \ell}(X^{[\ell]})$ is the
fundamental class of the closure of the following subvariety in
$X^{[\ell]}$:
\begin{eqnarray*}
\{ \xi_1 + \ldots + \xi_k|\,\, \xi_i \in L^{\nu_i}C_i \text{ for
every } i, \,\, \Supp(\xi_i) \cap \Supp(\xi_j) = \emptyset \text{
for } i \ne j \}.
\end{eqnarray*}
This follows from an induction on $k$ and $\ell(\nu_k)$, together
with (\ref{9.14.1}) and an argument similar to the proof of the
Theorem~9.14 in \cite{Na2}.
\end{remark}

\begin{remark}   \label{rmk_monodromy}
As noted by the referee, since the monodromy operators on the
cohomology of a surface preserve the Heisenberg operators on the
integral cohomology of Hilbert schemes, the operator $\mathfrak
m_{\la, \alpha}$ for a {\em non-algebraic} class $\alpha \in
H^2(X;\Z)/{\rm Tor}$ is integral if there exists a monodromy
sending $\alpha$ to a divisor. For instance, when $X$ is a K3
surface, it is known that for any class $\alpha \in H^2(X;\Z)$,
there exists a monodromy sending $\alpha$ to a divisor. Hence the
conclusions in Theorem~\ref{thm_pg_q} hold when $X$ is a K3
surface.
\end{remark}


%
%
%
%
%
%
%
\subsection{An algebraic model}
\label{sec:model}
\par
$\,$

The purpose of this subsection is to formalize the constructions
in the previous sections in a purely algebraic way.

Given a finite-dimensional graded Frobenius algebra $A$ over
$\mathbb Q$, one can construct a Fock space $\mathbb H_A$ of
Heisenberg algebra generated by $\mathfrak a_n(\alpha), n\in\Z,
\alpha \in A$. For the sake of simplicity, we assume that $A$ is
evenly graded, i.e., $A=\oplus_{i =0}^r A_{2i}$, with $A_0$ and
$A_{2r}$ being $1$-dimensional. Then, $\mathbb H_A$ is bigraded,
and its $n$-th component $A^{[n]}$ is still graded; in particular,
$A^{[1]}=A$. Furthermore, as graded vector spaces
$$\mathbb H_A \cong \mathbb H_{A_{\rm mid}}
\otimes \mathbb H_{A_0 \oplus A_{2r}}$$ where $A_{\rm mid}
=\oplus_{i =1}^{r-1} A_{2i}$ inherits a non-degenerate bilinear
form by restriction from $A$, $\mathbb H_{A_{\rm mid}}$ and
$\mathbb H_{A_0 \oplus A_{2r}}$ are the Fock spaces associated to
$A_{\rm mid}$ and $A_0 \oplus A_{2r}$ respectively. We denote by
$A_{\rm mid}^{[n]}$ the $n$-th component of $\mathbb H_{A_{\rm
mid}}$. Given an integral lattice $L_A$ in $A_{\rm mid}$ and an
integral basis $\alpha_1, \ldots, \alpha_k$ of it, we can define
the operators $\mathfrak a_{-\mu}(\alpha)$ associated to a
partition $\mu$ as before, and then define the operators
$\mathfrak m_{\la, \alpha}$ from the operators $\mathfrak
a_{-\mu}(\alpha)$ by declaring that the transition matrix is the
same as the (universal) one between the monomial symmetric
functions and the power-sum symmetric functions. Then we introduce
the elements $m_\la := \mathfrak m_{\la^1, \alpha_1} \ldots
\mathfrak m_{\la^k, \alpha_k}\vac \in \mathbb H_A$ associated to
$k$-tuple partitions $\la =(\la^1,\cdots,\la^k)$. We denote by
$L_A^{[n]}$ the lattice of $A_{\rm mid}^{[n]}$ which is the
$\Z$-span of the elements $m_\la$ with $|\la^1| + \cdots
|\la^k|=n$. Then by construction, the operators $\mathfrak
m_{\la^i, \alpha_i}, 1 \le i \le k,$ are ``integral", i.e. they
preserve $\oplus_{n \ge 0} L_A^{[n]}$. Furthermore, the
counterparts of Lemmas \ref{-a} and \ref{a=a1+a2} hold in the
current setup and they readily imply that the lattice $L_A^{[n]}$
is independent of the choice of the integral basis $\alpha_1,
\ldots, \alpha_k$ of $L_A$. If the determinant of the intersection
matrix of the lattice $L_A$ is $ \pm 1$, then the same arguments
as before imply that the determinant of the intersection matrix of
$L_A^{[n]}$ for each $n$ is $\pm 1$. However, it is not clear to
us whether or not the pairing on $L_A^{[n]}$ induced from the
integral pairing on $L_A$ is integral.

\begin{remark} \label{rmk_ref}
Assume further that we have an orthogonal direct sum $A_{\rm mid} =B
\oplus C$. Choose lattices $L_A, L_B$ and $L_C$ of $A_{\rm mid}$, $B$
and $C$ respectively such that $L_A =L_B \oplus L_C$.
As before, we can define the Fock spaces $\mathbb H_B$ and
$\mathbb H_C$, their $n$-th components $B^{[n]}$ and $C^{[n]}$,
and the lattices $L_B^{[n]}$ and $L_C^{[n]}$ for each $n$.
Clearly,
$\mathbb H_{A_{\rm mid}} \cong \mathbb H_{B} \otimes \mathbb H_{C}.$
Furthermore, by constructions the lattice structures are
compatible with the tensor product decomposition of Fock spaces:
\begin{eqnarray} \label{eq:tens}
\left( \oplus_{n \ge 0} L_A^{[n]} \right)
\cong \left( \oplus_{n \ge 0} L_B^{[n]}
\right) \bigotimes \left(\oplus_{n \ge 0} L_C^{[n]} \right).
\end{eqnarray}
\end{remark}
\vspace{.3cm}

On the other hand, take the lattice $M =\Z\cdot 1 \oplus \Z \cdot
x$, where $1$ is the unit of the Frobenius algebra $A$ and $x\in
A_{2r}$ has trace one. We can define an integral lattice $M^{[n]}$
in the $n$-th component of the Fock space $\mathbb H_{A_0 \oplus
A_{2r}}$ which is the $\Z$-span of $\frac1{\mathfrak z_\la}
\mathfrak a_{-\la}(1) \mathfrak a_{-\mu}(x)$, where $\la, \mu$ are
partitions such that $|\la| +|\mu| =n$. Now, $\big ( \oplus_{n \ge
0} L^{[n]} \big ) \otimes \big ( \oplus_{m \ge 0} M^{[m]} \big )$
is a lattice of the Fock space $\mathbb H_A =\mathbb H_{A_{\rm mid}}
\otimes \mathbb H_{A_0 \oplus A_{2r}}$. In other words,
$\oplus_{i=0}^n L^{[n-i]}  \otimes  M^{[i]}$ can be identified
with a lattice in $A^{[n]}$.

By taking $A =H^*(X)$, the Heisenberg algebra attains a geometric
meaning by Nakajima's construction \cite{Na2}. In this way, our
Theorem~\ref{thm_pg_q} can be restated that for a projective
surface $X$ with the given conditions, the lattice $\oplus_{i=0}^n
L^{[n-i]} \otimes M^{[i]}$ is identified with $H^*(\Xn; \mathbb
Z)/\text{\rm Tor}$. Our conjecture in the paragraph following
Theorem~\ref{thm_pg_q} can be further reformulated by
saying that $\oplus_{i=0}^n L^{[n-i]} \otimes  M^{[i]}$ can be
identified with $H^*(\Xn; \mathbb Z)/\text{\rm Tor}$ for every
simply-connected surface $X$.


%
%
%
%
%
\subsection{\bf Blown-up surfaces}
\label{subsect_e}
\par
$\,$

Let $\W X$ be the blown-up of $X$ at one point, and
$E$ be the exceptional curve. In this subsection,
we study a relation between integral bases
of $H^*(\W X^{[n]})$ and $H^*(X^{[n]})$.
Regard $H^*(X; \Z)/\text{Tor}$ as a sublattice of
$H^*(\W X; \Z)/\text{Tor}$, and $H^*(X)$ as a subspace
of $H^*(\W X)$.
Fix an integral basis $\{\alpha_1, \ldots, \alpha_k\}$ of
$\big ( H^1(X) \oplus H^2(X) \oplus H^3(X) \big )/\text{Tor}$.
Then, $\{\alpha_1, \ldots, \alpha_k, E\}$ is an integral basis
of $\big ( H^1(\W X) \oplus H^2(\W X) \oplus H^3(\W X)
\big )/\text{Tor}$.

Define ${\mathbb H}_{n, \W X}''$ to be the linear subspace of
${\mathbb H}_{n, \W X}' \subset \mathbb H_{n, \W X} = H^*(\W
X^{[n]})$ spanned by all Heisenberg monomial classes involving
only the creation operators $\mathfrak a_{-j}(\alpha)$ with
\begin{eqnarray}  \label{j>0}
j > 0, \qquad \alpha \,\, \in \,\, \bigoplus_{i=1}^3 H^i(X)
\,\, \subset \,\, \bigoplus_{i=1}^3 H^i(\W X).
\end{eqnarray}

Let $\W A \in H^*(\W X^{[n]})$. As in (\ref{A_la_mu}),
$\W A$ can be written as
\begin{eqnarray} \label{B_la_mu_nu}
\W A = \sum_{\la, \mu, \nu} \frac{1}{{\mathfrak z}_\la}
\cdot \mathfrak a_{-\la}(1)\mathfrak a_{-\mu}(x)
\mathfrak a_{-\nu}(E)(\W B_{\la, \mu, \nu})
\end{eqnarray}
where $\W B_{\la, \mu, \nu} \in {\mathbb H}_{n-|\la|-|\mu|-|\nu|,
\W X}''$ for partitions $\la$, $\mu$ and $\nu$.
Therefore, we see from (\ref{Phi_C2}) and (\ref{B_la_mu_nu}) that
$\W A$ can be further rewritten as
\begin{eqnarray} \label{A_la_mu_nu}
\W A = \sum_{\la, \mu, \nu} \frac{1}{{\mathfrak z}_\la}
\cdot \mathfrak a_{-\la}(1)\mathfrak a_{-\mu}(x)
\mathfrak m_{\nu, E}(\W A_{\la, \mu, \nu})
\end{eqnarray}
where $\W A_{\la, \mu, \nu} \in {\mathbb H}_{n-|\la|-|\mu|-|\nu|,
\W X}''$. The following is similar to Lemma~\ref{int_A_la_mu}.

\begin{proposition} \label{prop_b}
Let $\W A \in H^*(\W X^{[n]})$ be expressed as in
(\ref{A_la_mu_nu}). Then,
\begin{enumerate}
\item[{\rm (i)}] the classes $\W A_{\la, \mu, \nu}$ are
uniquely determined by $\W A$;

\item[{\rm (ii)}] $\W A$ is integral if and only if all
the classes $\W A_{\la, \mu, \nu}$ are integral.
\end{enumerate}
\end{proposition}
\begin{proof}
We follow an argument similar to the one used in the proof of
Lemma~\ref{int_A_la_mu}. First of all, since $1/\mathfrak z_\la
\cdot \mathfrak a_{-\la}(1), \mathfrak a_{-\mu}(x)$ and $\mathfrak
m_{\nu, E}$ are integral, $\W A$ is integral if all the classes
$\W A_{\la, \mu, \nu}$ are integral. So it remains to prove (i)
and the ``only if" part of (ii).

Next, in view of Lemma~\ref{int_A_la_mu},
it suffices to show that if
$\W A = \sum_{\nu} \mathfrak m_{\nu, E}(\W A_{\nu})$
where $\W A_{\nu} \in {\mathbb H}_{n-|\nu|, \W X}''$ for
partitions $\nu$, then we have
\par
$\quad${\rm (i$'$)} the classes $\W A_{\nu}$ are
uniquely determined by $\W A$;
\par
$\quad${\rm (ii$'$)} the classes $\W A_{\nu}$ are integral
when $\W A$ is integral.

\medskip
Let $n_0$ be the maximal integer such that $\W A_{\nu}
\ne 0$ for some $\nu \vdash n_0$.
Consider the action of all the
adjoint operators $(\mathfrak m_{\rho, E})^\dagger$,
$\rho \vdash n_0$ on $\W A$:
\begin{eqnarray} \label{prop_b.1}
(\mathfrak m_{\rho, E})^\dagger(\W A) = \sum_{\nu}
(\mathfrak m_{\rho, E})^\dagger \mathfrak m_{\nu, E}
(\W A_{\nu}), \qquad \rho \vdash n_0.
\end{eqnarray}
Recall from (\ref{Phi_C2}) that $\mathfrak m_{\rho, E}$
is a polynomial of the creation operators
$\mathfrak a_{-i}(E), \,\, i > 0$.
So $(\mathfrak m_{\rho, E})^\dagger$ is a polynomial of
the operators $\mathfrak a_i(E), \,\, i > 0$.
By the definition of ${\mathbb H}_{k, \W X}''$
and the fact that $(E, \alpha) = 0$ for $\alpha \in H^*(X)$,
if we write $\W A_{\nu} =
\mathfrak a_{\W A_{\nu}} \vac$ as in Proposition~\ref{prop_int},
then $\big [(\mathfrak m_{\rho, E})^\dagger,
\mathfrak a_{\W A_{\nu}} \big ] = 0$. Thus,
for $\rho \vdash n_0$, we see from (\ref{prop_b.1}) that
\begin{eqnarray} \label{prop_b.2}
& &(\mathfrak m_{\rho, E})^\dagger(\W A)
   = \sum_{\nu} \mathfrak a_{\W A_{\nu}}
   (\mathfrak m_{\rho, E})^\dagger \mathfrak m_{\nu, E}\vac
   = \sum_{|\nu| = n_0} \mathfrak a_{\W A_{\nu}}
   (\mathfrak m_{\rho, E})^\dagger \mathfrak m_{\nu, E}\vac
   \nonumber   \\
&=&\sum_{|\nu| = n_0} \big ( \mathfrak m_{\rho, E} \vac,
   \mathfrak m_{\nu, E}\vac \big )
   \cdot \mathfrak a_{\W A_{\nu}}\vac
   = \sum_{|\nu| = n_0} \big ( \mathfrak m_{\rho, E} \vac,
   \mathfrak m_{\nu, E}\vac \big ) \cdot \W A_{\nu}.
\end{eqnarray}

Applying Lemma~\ref{det=one} to $\underline{\beta} = \{ E \}$,
we see that the intersection matrix of the integral classes
$\mathfrak m_{\rho, E} \vac, \rho \vdash n_0$ is unimodular.
Therefore by (\ref{prop_b.2}), the classes
$\W A_{\nu}, \nu \vdash n_0$ are integral linear combinations
of the classes $(\mathfrak m_{\rho, E})^\dagger(\W A), \rho \vdash n_0$.
So the classes $\W A_{\nu}, \nu \vdash n_0$ are uniquely
determined by $\W A$. Moreover, in view of
Lemma~\ref{int_Heis}~(i), if the class $\W A$ is integral,
then all the classes
$(\mathfrak m_{\rho, E})^\dagger(\W A), \rho \vdash n_0$
are integral. Hence all the classes
$\W A_{\nu}, \nu \vdash n_0$ are integral as well.

Finally, repeating the above process to the class
\begin{eqnarray*}
\W A'
:= \W A - \sum_{|\nu| = n_0} \mathfrak m_{\nu, E}(\W A_{\nu})
= \sum_{|\nu| < n_0} \mathfrak m_{\nu, E}(\W A_{\nu}),
\end{eqnarray*}
we conclude that (i$'$) and (ii$'$) hold
for all the classes $\W A_{\nu}$.
\end{proof}

This proposition, together with (\ref{A_la_mu_nu}),
enables us to write down an
integral basis of $H^*(\W X^{[n]})$ whenever an integral basis of
the subspace ${\mathbb H}_{k, \W X}'' \cong {\mathbb H}_{k, X}'$
is known for every $k \le n$. It is also consistent with
subsection~\ref{sec:model}, since (\ref{eq:tens}) in
Remark~\ref{rmk_ref} is applicable to the orthogonal direct sum
$H^2(\W X; \Z) = H^2(X; \Z) \oplus \Z \cdot E$.

\end{document}